\documentclass[11pt,a4paper]{article}
\usepackage[T1]{fontenc}
\usepackage[utf8]{inputenc}
\usepackage[english]{babel}
\usepackage{lmodern,microtype}
\usepackage[margin=30mm]{geometry}
\usepackage{amsmath,amssymb,amsthm,mathtools,booktabs,array}
\usepackage{enumitem}
\usepackage{float}
\usepackage[hidelinks]{hyperref}
\hypersetup{pdftitle={On Semilocal Convergence of a Semismooth* Newton Method},
  pdfauthor={Tomáš Roubal}}
\theoremstyle{plain}
\newtheorem{theorem}{Theorem}[section]
\newtheorem{lemma}[theorem]{Lemma}
\newtheorem{proposition}[theorem]{Proposition}
\newtheorem{corollary}[theorem]{Corollary}
\theoremstyle{definition}
\newtheorem{definition}[theorem]{Definition}
\newtheorem{example}[theorem]{Example}
\theoremstyle{remark}
\newtheorem{remark}[theorem]{Remark}

\floatstyle{ruled}
\newfloat{algorithm}{htbp}{loa}[section]
\floatname{algorithm}{Algorithm}

\newcommand{\R}{\mathbb R}
\newcommand{\gph}{\operatorname{gph}}
\newcommand{\dist}{\operatorname{dist}}
\newcommand{\norm}[1]{\lVert#1\rVert}
\newcommand{\Areg}{\mathcal A_{\mathrm{reg}}F}

\newcommand{\bx}{\bar{x}}
\newcommand{\by}{\bar{y}}
\def\ball{{I\kern -.35em B}}

\begin{document}
\title{On Semilocal Convergence of a Semismooth$^*$ Newton Method}
\author{Tomáš Roubal\thanks{Corresponding author.
  Department of Computer Science, Faculty of Science, University of South
  Bohemia in České Budějovice, Na Zlaté stoce 1619, 370 05 České
  Budějovice, Czech Republic; Institute of Applied Mathematics,
  Graz University of Technology, Steyrergasse 30/III, 8010 Graz, Austria.
  Email: \href{mailto:roubalt@prf.jcu.cz}{\texttt{roubalt@prf.jcu.cz}}.}}
\date{}
\maketitle

\begin{abstract}
We give conditions under which the semismooth$^*$ Newton method solves
$0\in F(x)$ from a given initial point. At each step, we bound the
distance from the corrected point $(x_{k+1},0)$ to the graph of $F$
by a function of the correction length. We also bound
$\|A^{-1}B\|$ uniformly over the selected coderivative pairs $(A,B)$.
Together, these bounds give a test on the initial graph residual and
the radius of a prescribed ball. If the test holds, the iterates converge
to a solution in that ball. If the decrease modulus tends to zero, every
infinite run converges $q$-superlinearly.

For mixed complementarity problems on a box, we assume that the smooth
part has a Lipschitz continuous Jacobian and that every nonempty
principal submatrix is nonsingular. We bound the error after one step
for each regular pair. The bound is quadratic in the distance to a
solution and depends explicitly on the norm of the correction matrix.
For affine problems, we give a threshold on the norm of the vertical component of a
graph point. Below this threshold, all regular pairs at that point give
the same solution in one step, even when their rows come from different
branches. Hence any run with vanishing graph residuals terminates in
finitely many steps. For integer affine data, we give a test for an exact
step using only the problem data and the current iterate.
\end{abstract}

\noindent\textbf{Keywords:} semismooth$^*$ Newton method; semilocal convergence;
generalized equations; strong metric subregularity;
mixed complementarity problems.

\medskip
\noindent\textbf{Mathematics Subject Classification:} 49J53; 49M15; 65K15; 90C33.

\section{Introduction}
Gfrerer and Outrata~\cite{GO2021} introduced the semismooth$^*$ Newton
method for the inclusion
\begin{equation*}
  0\in F(x),
\end{equation*}
where $F:\R^n\rightrightarrows\R^n$ is a set-valued mapping with a closed graph. This model covers
variational inequalities and complementarity problems.
Gfrerer and Outrata proved local superlinear convergence when $F$ is
semismooth$^*$ at $(\bx,0)$ and strongly metrically regular around
$(\bx,0)$~\cite[Theorems~4.4 and~4.7]{GO2021}. Benko and
Mehlitz~\cite[Section~2.5]{BM} showed that strong metric subregularity
around the same point suffices under the same semismoothness$^*$ assumption.

We ask a semilocal question. Fix an initial point $x_0$ and a radius
$\rho>0$. When do the iterates converge to a solution in
$\ball[x_0,\rho]$? To answer this question, define the graph residual
\[
 r(x):=\dist((x,0),\gph F).
\]
At iteration $k$, choose a graph point $z_k=(\widehat x_k,\widehat y_k)$
and a regular coderivative pair $(A_k,B_k)$ at $z_k$. With
$v_k:=(x_{k+1}-\widehat x_k,-\widehat y_k)$, the correction gives
\[
 z_k+v_k=(x_{k+1},0).
\]
For the selected pairs and all sufficiently short corrections in the
prescribed neighborhood, we require
\[
 r(x_{k+1})
 =\dist(z_k+v_k,\gph F)
 \le \omega(\|v_k\|)\,\|v_k\|.
\]
The approximation step and a uniform bound on $\|A_k^{-1}B_k\|$ give
$\|v_k\|\le a r(x_k)$ with an explicit constant $a$. The residuals
therefore contract if $a\omega(a r(x_0))<1$. Theorem~\ref{thm:inner}
gives a radius bound and an a posteriori estimate of the distance to the
limit. If $\omega(s)\to0$ as $s\downarrow0$, convergence is
$q$-superlinear. Proposition~\ref{prop:uniform-inner-regularity} gives
a vanishing modulus for a suitable family of pairs under
semismoothness$^*$ and strong metric subregularity.

Gfrerer and Outrata~\cite{GO2022} also developed the SCD (subspace
containing derivative) variant and proved local superlinear convergence
under SCD regularity. Applications include Tresca friction and
variational inequalities~\cite{GOV2022Tresca,GOV2022}, as well as
Tikhonov functionals with nonsmooth nonconvex penalties~\cite{GHR2025}.
Gfrerer~\cite{G2025} gave a globally convergent method for nonsmooth
nonconvex optimization. Chen, Sun, and Zhang~\cite{CSZ2026} studied
two approximation steps, one based on quadratic programming and one on
a forward--backward step. They related the resulting methods to
G-semismooth Newton methods for suitable reformulations. Roubal and
Valdman~\cite{RV2025} extended the method to path following. Related
work treats polyhedral regularization~\cite{NVV2025}, total variation
regularization in tomography (Gfrerer et al.,
\href{https://arxiv.org/abs/2605.12041}{arXiv:2605.12041}), and
implicit programming (Gfrerer, Ko\v{c}vara, and Outrata,
\href{https://arxiv.org/abs/2412.05953}{arXiv:2412.05953}).

For nonlinear equations, the theorems of Kantorovich and
Smale~\cite{Kantorovich1948,Smale1986} give classical semilocal results.
Ferreira and Svaiter~\cite{FS2009,FS2012} used majorants to study exact
and inexact Newton methods. For Lipschitz continuous equations, Qi and
Sun~\cite[Theorem~3.3]{QS1993} proved a Kantorovich-type theorem for the
semismooth Newton method. Cibulka et al.~\cite{CDPRV} strengthened this
theorem and gave a corresponding result for the Josephy--Newton method.
Other approaches use Nachi--Penot derivatives~\cite{GYNP2022}, graphical
derivatives~\cite{HKMP2012}, or approximations of set-valued
mappings~\cite{DS2012,G2014,GG2013}.

Adly, Cibulka, and
Ngai~\cite{ACN2015} replace $\Psi$ in $0\in\Psi(x)+F(x)$ by a set-valued
approximation. Their proof uses the stability of metric regularity under
perturbations. Klatte and Kummer~\cite{KK2018} approximate $f$ in
$0\in f(x)+M(x)$. They use Newton maps to prove a Kantorovich-type
result under metric regularity. Az\'e and Chou~\cite{AC1995} start from
a small residual and construct graph points that converge to a solution.
They assume a surjective derivative with bounded inverse at a reference
point. In our method, $(A,B)$ defines an affine model of the whole graph
at $(\widehat x,\widehat y)$. The correction requires one linear solve
and allows singular $B$. We estimate how far the corrected point lies
from the graph.

For mixed complementarity problems on a box, the choice of pair matters.
Example~\ref{ex:box-pairs} gives a regular pair whose Newton step is zero
at a nonsolution arbitrarily close to a solution. For the affine problems
in Corollary~\ref{cor:box-fb}, all regular pairs at $(p,y)$ give the same
solution in one step whenever $\|y\|<\delta$. This includes pairs with
rows from different branches through $(p,y)$. Fischer and
Kanzow~\cite{FK1996} proved local one-step exactness for a Newton method
applied to the min-reformulation of a linear complementarity problem.
Our test uses $\|y\|$ and covers every pair in $\Areg(p,y)$.
For integer affine data, Corollary~\ref{cor:box-integral} gives a test
that needs only the problem data and the current iterate. It requires
neither a list of candidates nor a known solution.

We use Euclidean vector norms and the induced matrix norm. We write
$\norm{\cdot}_F$ for the Frobenius norm, $\langle\cdot,\cdot\rangle$ for
the Euclidean inner product, and $e_i$ for the $i$th standard basis
vector. The symbol $I_n$ denotes the identity matrix or the identity map
on $\R^n$. We write $\ball[x_0,\rho]$ and $\ball(x_0,\rho)$ for the
closed and open balls with center $x_0$ and radius $\rho$. Set
\[
 \dist(x,A):=\inf_{a\in A}\Vert x-a\Vert.
\]
We use $\inf\emptyset:=\infty$ and $\sup\emptyset:=0$.
For a set-valued mapping $F$, we write $\gph F$ for its graph.

Section~\ref{sec:prelim} recalls the semismooth$^*$ Newton method.
Section~\ref{sec:semilocal} proves the semilocal convergence theorem and
derives a vanishing modulus from regularity assumptions.
Section~\ref{sec:mcp} gives error bounds for box problems and a test
for one-step exactness in the affine case.

\section{Preliminaries and the Semismooth\texorpdfstring{$^*$}{*} Newton Method}\label{sec:prelim}
We use metric subregularity and strong metric subregularity at a
point~\cite{CDK2018,DR2014}. We also use the stronger property around a
point introduced in~\cite{GO2022}.

\begin{definition}\label{def:subreg}
Let $F:\R^n\rightrightarrows\R^n$ and $(\bx,\by)\in\gph F$.
\begin{enumerate}[label=\textup{(\roman*)}]
\item $F$ is \emph{metrically subregular} at $(\bx,\by)$ if some
$\kappa\ge0$ and a neighborhood $U$ of $\bx$ satisfy
\[
 \dist(x,F^{-1}(\by))\le\kappa\dist(\by,F(x))\quad\text{for each}\quad x\in U.
\]
\item $F$ is \emph{strongly metrically subregular} at $(\bx,\by)$ if
some $\kappa\ge0$ and a neighborhood $U$ of $\bx$ satisfy
\[
 \Vert x-\bx\Vert\le\kappa\dist(\by,F(x)) \quad\text{for each}\quad x \in U.
\]
\item\label{def:subreg-around} $F$ is \emph{strongly metrically subregular}
around $(\bx,\by)$ if there exist a neighborhood $W$ of $(\bx,\by)$ and
a constant $\kappa\ge0$ with the following property.
For each $(u,y)\in W\cap\gph F$, there is a neighborhood $U$ of $u$ such that
\[
 \Vert x-u\Vert\le\kappa\dist(y,F(x))\quad\text{for each}\quad x\in U.
\]
\end{enumerate}
\end{definition}

Fix a closed set $K\subset\R^n$ and a point $x\in K$. Following
\cite{Mordukhovich2006}, define the regular (Fréchet) normal cone by
\[
 \widehat N_K(x):=\Bigl\{x^*\in\mathbb{R}^n:\
 \limsup_{v\to x,\ v\in K\setminus\{x\}}
 \frac{\langle x^*,v-x\rangle}{\Vert v-x\Vert}\le0\Bigr\},
\]
and the limiting (Mordukhovich) normal cone by
\[
 N_K(x):=\bigl\{x^*\in\mathbb{R}^n:\ \exists\,x_k\in K,\ x_k^*\in\widehat N_K(x_k)
 \text{ with } x_k\to x,\ x_k^*\to x^*\bigr\}.
\]
For $x\notin K$, set $\widehat N_K(x)=N_K(x)=\emptyset$. If $K$ is
convex and $x\in K$, the limiting normal cone takes the form
\[
 N_K(x)=\{x^*\in\R^n:\langle x^*,v-x\rangle\le0
                  \text{ for all }v\in K\}.
\]

\begin{definition}
Let $F:\R^n\rightrightarrows\R^n$ have a closed graph and take
$(\bx,\by)\in\gph F$. For $v\in\R^n$, the \emph{limiting coderivative}
of $F$ at $(\bx,\by)$ is
\[
 D^*F(\bx,\by)(v):=\{u\in\R^n:(u,-v)\in N_{\gph F}(\bx,\by)\}.
\]
\end{definition}

For $(p,y)\in\gph F$, the set $\Areg(p,y)$ consists of all pairs
$(A,B)\in\R^{n\times n}\times\R^{n\times n}$ such that $A$ is
nonsingular and
\begin{equation*}
 A^\top e_i\in D^*F(p,y)(B^\top e_i),\quad i=1,\ldots,n.
\end{equation*}
These pairs define the Newton step in the following version of
\cite[Algorithm~2]{GO2021}.

\begin{algorithm}[H]
\caption{Semismooth$^*$ Newton method}
\label{alg:Semismooth-star}
Choose $x_0\in\R^n$ and $\beta\ge1$. For $k=0,1,\ldots$, take the
following steps.
\begin{enumerate}[leftmargin=1.7em]
\item If $0\in F(x_k)$, stop.
\item Choose $(\widehat x_k,\widehat y_k)\in\gph F$ such that
$\Areg(\widehat x_k,\widehat y_k)\ne\emptyset$ and
\[
 \Vert(\widehat x_k-x_k,\widehat y_k)\Vert\le
 \beta\dist((x_k,0),\gph F).
\]
\item Choose $(A_k,B_k)\in\Areg(\widehat x_k,\widehat y_k)$ and solve
$A_kd_k=-B_k\widehat y_k$.
\item Set $x_{k+1}:=\widehat x_k+d_k$.
\end{enumerate}
\end{algorithm}

\begin{definition}
Let $F:\R^n\rightrightarrows\R^n$ and $(\bx,\by)\in\gph F$.
We call $F$ \emph{semismooth}$^*$ at $(\bx,\by)$ if for every
$\varepsilon>0$ there is a neighborhood $W$ of $(\bx,\by)$ such that
\[
 |\langle x^*,x-\bx\rangle-\langle y^*,y-\by\rangle|
 \le\varepsilon\Vert(x^*,y^*)\Vert\,\Vert(x,y)-(\bx,\by)\Vert
\]
for all $(x,y)\in\gph F\cap W$ and $(y^*,x^*)\in\gph D^*F(x,y)$.
\end{definition}
Gfrerer and Outrata~\cite[Corollary~3.3]{GO2021} show that this condition
is equivalent to their original definition.

\section{A Semilocal Convergence Theorem}\label{sec:semilocal}
Let $F:\R^n\rightrightarrows\R^n$ have a nonempty closed graph $G:=\gph F$.
We seek a solution $\bx\in\R^n$ of
\begin{eqnarray}
\label{eq:problem-semilocal}
    0\in F(\bx).
\end{eqnarray}
Define the graph residual by
\begin{equation*}
 r(x):=\dist((x,0),G) \quad\text{for} \quad x\in\R^n.
\end{equation*}
The function $r$ is Lipschitz on $\R^n$.
Since $G$ is closed, $r(x)=0$ exactly when $0\in F(x)$.
\begin{lemma}\label{lem:residual}
Let $x_0\in\R^n$, $\rho>0$, $a>0$, $\gamma\in [0,1)$, and $r_0:=r(x_0)>0$.
Let $\phi:[0,r_0]\longrightarrow[0,r_0]$ satisfy $\phi(0)=0$ and
\[
 \phi(t)\le\gamma\, t\quad\text{for each}\quad  0\le t\le r_0.
\]
Suppose the iterates starting at $x_0$ satisfy
\begin{equation*}
 \norm{x_{k+1}-x_k}\leq ar(x_k)\quad\text{and}\quad r(x_{k+1})\leq\phi(r(x_{k}))\quad\text{for}\quad k=0,1,2,\dots
\end{equation*}
If
\begin{equation*}
 \frac{ar_0}{1-\gamma}\le\rho,
\end{equation*}
then all iterates lie in $\ball[x_0,\rho]$ and converge to a solution
$\bx$ of \eqref{eq:problem-semilocal}. For each $k\ge0$,
\begin{equation}\label{eq:tail}
 r(x_k)\leq\gamma^k r_0\quad\text{and}\quad
 \norm{x_k-\bx}\le\frac{a}{1-\gamma}r(x_k).
\end{equation}
If $\phi(t)/t\to0$ as $t\downarrow0$, convergence is $q$-superlinear.
If $\phi(t)=\lambda t^{1+\alpha}$ with $\lambda,\alpha>0$ and
$\lambda r_0^\alpha<1$, take $\gamma=\lambda r_0^\alpha$. Then
\begin{equation}\label{eq:qpower}
 \norm{x_{k+1}-\bx}
 \le\frac{a\lambda}{1-\gamma}\norm{x_k-\bx}^{1+\alpha}.
\end{equation}
\end{lemma}
\begin{proof}
The residual bound gives $r(x_k)\le\gamma^k r_0$ by induction.
Summing the step lengths gives
\[
 \norm{x_{k+1}-x_0}
 \le a\sum_{j=0}^k r(x_j)
 \le\frac{ar_0}{1-\gamma}\le\rho.
\]
Thus all iterates lie in $\ball[x_0,\rho]$. For integers $k>s\ge0$,
\begin{eqnarray*}
 \norm{x_{k}-x_s}
 &\leq&  \norm{x_{k}-x_{k-1}} + \cdots +\norm{x_{s+1}-x_s}\leq  a r(x_{k-1})+\cdots + a r(x_{s})\\
 &\leq& a\frac{1-\gamma^{k-s}}{1-\gamma}r(x_s)\leq  \frac{a r(x_s)}{1-\gamma}.
 \end{eqnarray*}
As $s\to\infty$, $r(x_s)\le\gamma^s r_0\to0$. Hence $(x_k)$ is Cauchy.
Its limit $\bx$ lies in $\ball[x_0,\rho]$ and satisfies $r(\bx)=0$
by continuity. Let $k\to\infty$ with $s$ fixed to get~\eqref{eq:tail}.
Apply the error bound at $k+1$ to obtain
\[
 \norm{x_{k+1}-\bx}
 \le\frac{a}{1-\gamma}\phi(r(x_k)).
\]
Together with $r(x_k)\le\norm{x_k-\bx}$, this gives both the
$q$-superlinear rate and~\eqref{eq:qpower} under the stated assumptions.
\end{proof}

Estimate~\eqref{eq:tail} gives $r$-linear convergence. The residual
factor $\gamma$ need not bound the ratio of successive errors.
The ball may contain several solutions.

Fix $x_0\in\R^n$, $\rho>0$, $\beta\ge1$, and $c\ge0$.
Let $r_0:=r(x_0)>0$ and set
\begin{equation*}
 \mathcal Z:=G\cap\ball[(x_0,0),\rho+\beta r_0]\quad\text{and}
 \quad a:=\beta\sqrt{1+c^2}.
\end{equation*}
For each $z\in\mathcal Z$, choose a nonempty set
$\mathcal R(z)\subseteq\Areg(z)$ such that
\begin{equation}\label{eq:Cbound}
 \norm{A^{-1}B}\leq c
 \quad\text{for all}\quad(A,B)\in\mathcal R(z).
\end{equation}

At each step of Algorithm~\ref{alg:Semismooth-star} with
$x\in\ball[x_0,\rho]$ and $r(x)\le r_0$, choose a graph point
$z=(\widehat x,\widehat y)\in\mathcal Z$ and a pair $(A,B)$ satisfying
\begin{equation}\label{eq:projection}
 \norm{z-(x,0)}\leq\beta r(x)
 \quad\text{and}\quad (A,B)\in\mathcal R(z).
\end{equation}
The closed graph $G$ contains a nearest point to $(x,0)$. This point
meets the distance bound in~\eqref{eq:projection} even for $\beta=1$.
Every graph point $z$ that meets the bound lies in $\mathcal Z$, since
\[
 \norm{z-(x_0,0)}
 \leq \norm{z-(x,0)}+\norm{x-x_0}
 \leq \beta r_0+\rho.
\]
We can therefore choose a pair from $\mathcal R(z)$ at each such point.

\begin{theorem}\label{thm:inner}
 Let $\omega:[0,ar_0]\longrightarrow[0,\infty)$
be nondecreasing. For every $z=(p,y)\in\mathcal Z$ and every
$(A,B)\in\mathcal R(z)$, set $v:=(-A^{-1}By,-y)$ and suppose that
\begin{equation}\label{eq:inner}
 \operatorname{dist}(z+v,G)
 \le \omega(\lVert v\rVert)\lVert v\rVert
 \quad\text{whenever}\quad\lVert v\rVert\le ar_0.
\end{equation}
If
\begin{equation}\label{eq:innerK}
 \gamma:=a\omega(ar_0)<1\quad\text{and}
 \quad \frac{ar_0}{1-\gamma}\le\rho,
\end{equation}
then every run of Algorithm~\ref{alg:Semismooth-star} from $x_0$ using
\eqref{eq:projection} stays in $\ball[x_0,\rho]$. It either reaches a
solution $\bx$ in finitely many steps or converges to a solution $\bx$
in $\ball[x_0,\rho]$. In both cases, \eqref{eq:tail} bounds the residual and
the error at each iterate. If $\omega(s)\to0$ as $s\downarrow0$, every
infinite run converges $q$-superlinearly.
\end{theorem}

\begin{proof}
Suppose $x_k\in\ball[x_0,\rho]$ and $0<r(x_k)\le r_0$. Choose
$z_k=(\widehat x_k,\widehat y_k)\in\mathcal Z$ and $(A_k,B_k)$ by
\eqref{eq:projection}. Put $C_k:=A_k^{-1}B_k$. The Newton update
$x_{k+1}=\widehat x_k-C_k\widehat y_k$ gives
\[
 \begin{aligned}
\norm{x_{k+1}-x_k}\leq \norm{\widehat{x}_k-x_k}+c\norm{\widehat{y}_k}
 \leq \sqrt{1+c^2}\,
       \norm{(\widehat{x}_k-x_k,\widehat{y}_k)}
 \leq a r(x_k).
 \end{aligned}
\]
For $v_k:=(-C_k\widehat y_k,-\widehat y_k)$, we have
\[
 z_k+v_k=(x_{k+1},0)
 \quad\text{and}\quad
 \norm{v_k}\leq\sqrt{1+c^2}\norm{\widehat{y}_k}
 \leq ar(x_k)\le ar_0.
\]
Apply \eqref{eq:inner} and use the monotonicity of $\omega$ to get
\[
 r(x_{k+1})
 =\dist(z_k+v_k,G)
 \leq\omega(\norm{v_k})\norm{v_k}
 \le ar(x_k)\omega(ar(x_k)).
\]
Set $\phi(t):=at\omega(at)$. This function is nondecreasing on $[0,r_0]$,
satisfies $\phi(0)=0$, and obeys
\[
 \phi(t)\le a\omega(ar_0)t=\gamma t
 \quad\text{for each}\quad 0\le t\le r_0.
\]
We now prove $x_k\in\ball[x_0,\rho]$ and $r(x_k)\le\gamma^k r_0$ by
induction. Both claims hold at $k=0$. Suppose they hold through $k$.
If the algorithm takes another step, the preceding estimates give
$r(x_{k+1})\le\gamma r(x_k)\le\gamma^{k+1}r_0$ and
\[
 \norm{x_{k+1}-x_0}
 \le a\sum_{j=0}^k r(x_j)
 \le a r_0\sum_{j=0}^k\gamma^j
 \le \frac{a r_0}{1-\gamma}\le\rho.
\]
If the algorithm stops, repeat its final iterate. The residual is
then zero, so the step and residual bounds still hold. Apply
Lemma~\ref{lem:residual} to get convergence and~\eqref{eq:tail}.
If $\omega(s)\to0$ as $s\downarrow0$, then
$\phi(t)/t=a\omega(at)\to0$ as $t\downarrow0$.
The lemma gives $q$-superlinear convergence for every infinite run.
\end{proof}

We call \eqref{eq:inner} the \emph{decrease condition} and $\omega$
its \emph{modulus}.

\begin{corollary}\label{cor:innerpower}
In Theorem~\ref{thm:inner}, either of the following bounds can replace
\eqref{eq:inner}.
\begin{enumerate}[label=\textup{(\roman*)},leftmargin=2em]
\item If $\dist(z+v,G)\le\varepsilon\norm v$,
      $a\varepsilon<1$, and $ar_0/(1-a\varepsilon)\le\rho$, then the convergence is $r$-linear.
\item If, for $0<\alpha\le1$,
      $\dist(z+v,G)\le M\norm v^{1+\alpha}$, set
      $\lambda:=Ma^{1+\alpha}$. The conditions
      \[
       \gamma:=\lambda r_0^\alpha<1\quad\text{and}\quad
       \frac{ar_0}{1-\gamma}\le\rho
      \]
      ensure $q$-order at least $1+\alpha$ and~\eqref{eq:qpower}.
\end{enumerate}
\end{corollary}

\begin{proof}
Substitute $\omega(s)=\varepsilon$ and $\omega(s)=Ms^\alpha$, respectively,
into Theorem~\ref{thm:inner} and Lemma~\ref{lem:residual}.
\end{proof}

An upper bound $\bar r_0\ge r(x_0)$ suffices in Theorem~\ref{thm:inner}. Use $\bar r_0$ in place
of $r_0$ when defining $\mathcal Z$ and checking \eqref{eq:inner} and
\eqref{eq:innerK}. This enlarges the region and the range of step lengths
to check. On the enlarged set $\mathcal Z$, the families $\mathcal R(z)$ must
also be nonempty and satisfy \eqref{eq:Cbound} with the chosen $c$. Lemma~\ref{lem:residual} and Theorem~\ref{thm:inner} then
give~\eqref{eq:tail} with $\bar r_0$ in place of $r_0$. For instance,
any $(\widehat x,\widehat y)\in G$ gives the bound
$\bar r_0:=\Vert(\widehat x-x_0,\widehat y)\Vert$.
Example~\ref{ex:cubic} uses only $r_0\le\varepsilon$.

We now derive \eqref{eq:Cbound} and \eqref{eq:inner} from regularity
of $F$. Theorem~\ref{thm:inner} then applies whenever \eqref{eq:innerK}
holds.

\begin{proposition}\label{prop:uniform-inner-regularity}
Suppose that $F$ is strongly metrically subregular around every point of
$\mathcal Z$ and semismooth$^*$ at every $(\bx,0)\in\mathcal Z$.
There exists $c\ge0$ such that the sets
\[
 \mathcal R(z):=
 \left\{
 (A,B)\in\mathcal A_{\mathrm{reg}}F(z):
 \begin{array}{l}
 \Vert A^{-1}B\Vert \le c,\\
 \Vert A^{-1}\Vert \Vert [A\ B]\Vert _F\le\sqrt{n(1+c^2)}
 \end{array}
 \right\},
 \quad z\in\mathcal Z,
\]
are nonempty. For every $s_0>0$, there is a nondecreasing function
$\omega:[0,s_0]\to[0,1]$ with
\[
 \omega(0)=0\quad\text{and}\quad
 \lim_{s\downarrow0}\omega(s)=0,
\]
such that
\[
 \operatorname{dist}((p-A^{-1}By,0),G)
 \le
 \omega(\Vert (A^{-1}By,y)\Vert )\,\Vert (A^{-1}By,y)\Vert
\]
for all $(p,y)\in\mathcal Z$ and all $(A,B)\in\mathcal R(p,y)$
with $\Vert(A^{-1}By,y)\Vert\le s_0$.
Taking $s_0=ar_0$ gives \eqref{eq:Cbound} and \eqref{eq:inner}.
\end{proposition}

\begin{proof}
By compactness, finitely many outer neighborhoods in
Definition~\ref{def:subreg}\ref{def:subreg-around} cover $\mathcal Z$.
Choose $\kappa_0>0$ at least as large as their subregularity moduli.
The same bound then holds throughout $\mathcal Z$. The inner
neighborhood of $p$ may depend on $(p,q)$. See~\cite{Roubal2025}
for related uniformity results. At each $z\in\mathcal Z$, the proof
of \cite[Proposition~2.8]{BM} gives a pair
$(I_n,B)\in\mathcal A_{\mathrm{reg}}F(z)$ such that each row of $B$ has
norm at most $\kappa_0$. Hence
\[
 \Vert B\Vert \le\Vert B\Vert _F\le\sqrt n\,\kappa_0.
\]
Set $c:=\sqrt n\,\kappa_0$. Then $(I_n,B)\in\mathcal R(z)$ because
\[
 \Vert [I_n\ B]\Vert _F^2
 =n+\Vert B\Vert _F^2
 \le n(1+\Vert B\Vert ^2)
 \le n(1+c^2).
\]

Fix $s_0>0$. For $0\le t\le s_0$, set
\[
 \omega(t):=
 \sup\left\{
 \frac{r(p-A^{-1}By)}
      {\Vert (A^{-1}By,y)\Vert }:
 \begin{array}{l}
 (p,y)\in\mathcal Z,\quad
 (A,B)\in\mathcal R(p,y),\\
 0<\Vert (A^{-1}By,y)\Vert \le t
 \end{array}
 \right\}.
\]
The point $(p,y)\in G$ gives the bound
$r(p-A^{-1}By)\le\Vert(A^{-1}By,y)\Vert$.
Thus $0\le\omega(t)\le 1$ for $0\leq t \leq s_0$, $\omega(0)=0$, and $\omega$ is nondecreasing.
Its definition gives the required residual estimate.

Suppose $\omega(t)\not\to0$ as $t\downarrow0$. Then there are
$\varepsilon>0$, points $(p_k,y_k)\in\mathcal Z$, and pairs
$(A_k,B_k)\in\mathcal R(p_k,y_k)$ such that
\[
 s_k:=\Vert (A_k^{-1}B_ky_k,y_k)\Vert >0,
 \quad s_k\to0,
 \quad\text{and}\quad
 r(p_k-A_k^{-1}B_k y_k)\ge\varepsilon s_k.
\]
Since $\mathcal Z$ is compact and $\Vert y_k\Vert\le s_k\to0$, a
subsequence of $(p_k,y_k)$ converges to $(\bar p,0)\in\mathcal Z$. Keep the same
indices for this subsequence. Then $\bar p$ solves
\eqref{eq:problem-semilocal}, and $F$ is semismooth$^*$ at $(\bar p,0)$.
Strong metric subregularity at $(\bar p,0)$ gives, for all large $k$,
\[
 d_k:=\Vert (p_k-\bar p,y_k)\Vert
 \le\sqrt{1+\kappa_0^2}\,\Vert y_k\Vert
 \le\sqrt{1+\kappa_0^2}\,s_k.
\]
Set
\[
 \eta:=\frac{\varepsilon}
 {2\sqrt{n(1+c^2)(1+\kappa_0^2)}}.
\]
By semismoothness$^*$ at $(\bar p,0)$ and
\cite[Proposition~4.3]{GO2021}, for all large $k$ we have
\[
 \begin{aligned}
 \dist((p_k-A_k^{-1}B_ky_k,0),G)
 &\le\Vert p_k-\bar p-A_k^{-1}B_ky_k\Vert  \le\eta\Vert A_k^{-1}\Vert \Vert [A_k B_k]\Vert _F\,d_k\\
 &\le\eta\sqrt{n(1+c^2)(1+\kappa_0^2)}\,s_k
 =\frac{\varepsilon}{2}s_k.
 \end{aligned}
\]
This contradicts the lower bound $\varepsilon s_k$.
\end{proof}

\begin{remark}\label{rem:global-semismooth-star}
If $\gph F$ is a finite union of closed convex sets or a closed
subanalytic set, then $F$ is semismooth$^*$ at each point of its graph.
See \cite[Proposition~3.5]{GO2021} and \cite[Proposition~2.10]{GO2022}.
Either condition gives the semismoothness$^*$ assumption in
Proposition~\ref{prop:uniform-inner-regularity}.
\end{remark}

\begin{lemma}\label{lem:box-subregularity}
Let $f:\R^n\longrightarrow\R^n$ be continuously differentiable and let
$K:=\prod_{i=1}^n[\ell_i,u_i]$, where $\ell_i<u_i$ for each $i$.
Suppose that every nonempty principal submatrix of $\nabla f(p)$
is nonsingular for every $p\in K$. Then $F:=f+N_K$ is strongly
metrically subregular around every point of its graph and
semismooth$^*$ at every point of its graph.
\end{lemma}

\begin{proof}
Fix $(p,y)\in\gph F$ and let
\[
 c(p):=\max_{\varnothing\ne J\subseteq\{1,\ldots,n\}}
          \|(\nabla f(p)_{JJ})^{-1}\|.
\]
Let $(h,k)$ be a tangent vector to $\gph F$ at $(p,y)$.
Since the graph is a finite union of closed branches, a sequence
defining $(h,k)$ has a subsequence on one branch through $(p,y)$.
Let $J$ consist of the indices for which the chosen branch imposes
$y_i=f_i(p)$. For each $i\notin J$, the branch fixes $p_i$ at one
of its bounds. Hence
\[
h_{J^c}=0\quad\text{and}\quad k_J=\nabla f(p)_{JJ}h_J.
\]
If $J=\varnothing$, then $h=0$. Otherwise, inversion of the
principal block gives $\|h\|\le c(p)\|k\|$.

By the graphical derivative criterion
\cite[Theorem~4E.1]{DR2014}, every $\kappa>c(p)$ is a modulus of
strong metric subregularity at $(p,y)$. Now fix
$(\bar p,\bar y)\in\gph F$ and $\kappa>c(\bar p)$. Since $c$ is
continuous on $K$, the same $\kappa$ is admissible at every graph
point $(p,y)$ sufficiently close to $(\bar p,\bar y)$. This proves
strong metric subregularity around $(\bar p,\bar y)$.

Finally, $\gph N_K$ is a finite union of convex polyhedra. Thus
$N_K$ is semismooth$^*$ at every graph point by
Remark~\ref{rem:global-semismooth-star}. The same holds for
$f+N_K$ by \cite[Proposition~3.6]{GO2021}.
\end{proof}

Lemma~\ref{lem:box-subregularity} gives the regularity assumptions
of Proposition~\ref{prop:uniform-inner-regularity} for box
problems. Under the additional hypotheses of
Corollary~\ref{cor:box-algorithm}\ref{it:box-semilocal}, the
decrease condition~\eqref{eq:inner} holds with $\omega(s)=2c^3Ls$.
\begin{example}\label{ex:cubic}
Take $\varepsilon:=1/100$ and consider $0\in F(x)$, where
\[
 f_\sigma(x):=\sigma(x+x^3)-\varepsilon\quad\text{and}\quad
 F(x):=\{f_{-1}(x),f_1(x)\}\quad \text{for}\quad \sigma\in\{-1,1\}.
\]
Set $x_0:=0$, $\rho:=1/4$, $\beta:=\sqrt5$, and $c:=1$.
Then $a=\sqrt{10}$. Since $G$ is closed, $(0,-\varepsilon)\in G$,
and $(0,0)\notin G$, we have $0<r_0\le\varepsilon$.

For $(p,y)\in G$, set
\[
 \mathcal R(p,y):=
 \{(f_\sigma'(p),1):\sigma\in\{-1,1\},\ y=f_\sigma(p)\}.
\]
The branches meet only at $(0,-\varepsilon)$. At each other point on a branch,
$(f_\sigma'(p),-1)$ is a regular normal to $G$. Taking limits along
that branch gives a limiting normal at the intersection. Thus
$\mathcal R(p,y)\subset\Areg(p,y)$. Since
$|f_\sigma'(p)|=1+3p^2\ge1$, \eqref{eq:Cbound} holds with $c=1$.

For the approximation step, set $\widehat x:=x$ and choose
$\widehat y\in F(x)$ with the smallest absolute value.
Suppose $|x|\le\rho$ and $r(x)\le r_0$. Take a nearest point
$(p,y)\in G$ to $(x,0)$. Then $|p|\le\rho+r_0<1/2$.
Both branches are Lipschitz with constant $2$ on $[-1/2,1/2]$.
Choose $\sigma$ with $y=f_\sigma(p)$. Thus
\[
 \norm{(\widehat x-x,\widehat y)}
 \leq |f_\sigma(x)|
 \leq |y|+2|p-x|
 \leq\sqrt{5}\,r(x).
\]
Hence $(\widehat x,\widehat y)$ satisfies \eqref{eq:projection}.

To check \eqref{eq:inner}, take $z=(p,y)\in\mathcal Z$ and
$(A,B)=(f_\sigma'(p),1)\in\mathcal R(p,y)$. Put $p^+:=p-y/f_\sigma'(p)$
and $v:=(p^+-p,-y)$. If $\norm v\le ar_0$, then $p,p^+\in[-1,1]$ because
\[
 \max\{|p|,|p^+|\}\le \rho+(\beta+a)r_0<1.
\]
Taylor's estimate, $f_\sigma(p)+f_\sigma'(p)(p^+-p)=0$, and
$|f_\sigma''|\le6$ on $[-1,1]$ give
\[
 \dist(z+v,G)
 \le |f_\sigma(p^+)|
 \le 3|p^+-p|^2
 \le 3\norm v^2.
\]
Thus $\omega(s):=3s$ satisfies \eqref{eq:inner}. We also have
\[
 \gamma=3a^2r_0\le\frac{3}{10}
 \quad \text{and}\quad
 \frac{ar_0}{1-\gamma}\le\frac{\sqrt{10}}{70}<\rho.
\]
These inequalities give \eqref{eq:innerK}. Corollary~\ref{cor:innerpower}(ii) shows
that the method either stops at a solution or converges to one with
$q$-order at least two.
\end{example}

\section{Mixed Complementarity Problems on a Box}\label{sec:mcp}
We apply Theorem~\ref{thm:inner} to the inclusion
\begin{equation}\label{eq:box-mcp}
 0\in F(x):=f(x)+N_K(x),\quad K:=\prod_{i=1}^n[\ell_i,u_i],
\end{equation}
where $\ell_i<u_i$ for each $i$. As in Section~\ref{sec:semilocal},
set $G:=\gph F$ and $r(x):=\dist((x,0),G)$. Assume that $G$ is closed
and that $f$ is continuously differentiable on an open set $O$.
Assume that every nonempty principal submatrix of $\nabla f(p)$ is
nonsingular for every $p\in O$.

Choose a pattern $\sigma\in\{-1,0,+1\}^n$, where $-1$, $0$, and $+1$
denote the lower-bound, free, and upper-bound cases, respectively.
A point $(p,y)$ belongs to the branch $P_\sigma$ if, for each $i$,
it satisfies the conditions in the middle column of the row indexed
by $\sigma_i$. A point $h\in O$ is a \emph{candidate} for $\sigma$
if it satisfies the corresponding equation in the last column for every $i$.
\[
\begin{array}{c@{\quad}|@{\quad}l@{\quad}|@{\quad}l}
 \sigma_i & (p,y)\in P_\sigma & \text{candidate equation}\\ \hline
 -1 & p_i=\ell_i,\quad y_i\le f_i(p)         & h_i=\ell_i\\
  0 & \ell_i\le p_i\le u_i,\quad y_i=f_i(p) & f_i(h)=0\\
 +1 & p_i=u_i,\quad y_i\ge f_i(p)           & h_i=u_i
\end{array}
\]
If $p_i$ is at a bound and $y_i=f_i(p)$, both the free label and
the label for that bound are admissible for coordinate $i$.

For the chosen pattern, let $J:=\{i:\sigma_i=0\}$ denote the set of
free indices. For $p\in O$ and $J\ne\varnothing$, define
$C_J(p)\in\R^{n\times n}$ to have $(\nabla f(p)_{JJ})^{-1}$ as its
$J\times J$ block and zeros in all other entries. Set
$C_\varnothing(p):=0$.

In the affine case $f(x)=Mx+q$, take $O=\R^n$ and write $C_J$ for
the constant matrix $C_J(p)$. Each pattern then has a unique
candidate $h^\sigma$.

Fix a closed convex set $D\subset O$. Denote the candidates and
solutions in $D$ by
\[
 \mathcal H:=\{h\in D:h\text{ is a candidate for some }\sigma\}\quad\text{and}
 \quad \mathcal S:=D\cap F^{-1}(0).
\]
Every solution solves the equations for at least one pattern, so
$\mathcal S\subseteq\mathcal H$. Let $\Pi_K$ denote the Euclidean projection onto $K$. For $t>0$, set
\[
 T(x):=\Pi_K(x-tf(x))\quad \text{and}\quad \eta(x):=\|x-T(x)\|.
\]
Then $\mathcal S=\{h\in\mathcal H:\eta(h)=0\}$.

Fix $c>0$ and $L\ge0$. In Proposition~\ref{prop:box-exact}, $c$
bounds $\|C_J(p)\|$ for all $p\in D$ and all $J$, and $L$ is
a Lipschitz constant for $\nabla f$ on $D$. Set $m:=c$ if $L=0$ and $m:=2c$ otherwise.
Define
\[
 L_T:=\sup_{x\in D}\|I_n-t\nabla f(x)\|\quad\text{and}\quad b:=(1+L_T)m+t.
\]
With $\min\varnothing:=\infty$ and $1/0:=\infty$, set
\[
 \delta:=\min\left\{
 \min_{h\in\mathcal H\setminus\mathcal S}\frac{\eta(h)}b,\quad
 \min_{\substack{h,h'\in\mathcal S\\h\ne h'}}\frac{\|h-h'\|}{2m},\quad
 \frac{1}{4c^2L}\right\}.
\]

\begin{proposition}\label{prop:box-exact}
Assume that $\|C_J(p)\|\le c$ for all $p\in D$ and all $J$.
Let $\nabla f$ be Lipschitz continuous on $D$ with constant $L$. If $L>0$,
assume also that $D$ is bounded. Then $\mathcal H$ is finite and
$\delta>0$. The following estimates hold.
\begin{enumerate}[label=(\alph*)]
\item\label{it:box-pair}
If $(p,y)\in P_\sigma$ and $p\in O$, then $(I_n,C_J(p))\in\Areg(p,y)$.
\item\label{it:box-exact}
Suppose $(p,y)\in G$, $\|y\|<\delta$, and $\ball[p,m\|y\|]\subset D$.
All branches through $(p,y)$ contain a common point $(h,0)$ with
$h\in\mathcal S$ and $\|p-h\|\le m\|y\|$. For any $(A,B)\in\Areg(p,y)$,
the Newton step $x^+:=p-A^{-1}By$ satisfies
\begin{equation}\label{eq:box-exact-correction}
 \|x^+-h\|\le \tfrac12 L\|A^{-1}B\|\,\|p-h\|^2.
\end{equation}
\end{enumerate}
\end{proposition}

\begin{proof}
Fix a pattern. The Jacobian of its candidate equations has rows
$e_i^\top$ for $i\notin J$ and $\nabla f_i(p)$ for $i\in J$.
The principal submatrix assumption makes this Jacobian nonsingular.
Each candidate is therefore isolated by the inverse function theorem. If $D$ is
bounded, the candidates for this pattern form a closed subset of the
compact set $D$, so there are finitely many. If $D$ is unbounded,
$L=0$ and $f$ is affine on $D$. Each pattern then has at most one
candidate. Thus $\mathcal H$ is finite. Since $\eta>0$ on
$\mathcal H\setminus\mathcal S$, we get $\delta>0$.

For part~\ref{it:box-pair}, first take $(q,z)\in P_\sigma$ with
$q\in O$ and each scalar pair in the relative interior of its piece.
Near $(q,z)$ the graph is a $C^1$ manifold. Each tangent vector
$(a,b)$ satisfies $a_i=0$ for $i\notin J$ and
$b_J=\nabla f(q)_{JJ}a_J$. Thus $a=C_J(q)b$.
It follows that $(u,-C_J(q)^\top u)$ is a regular normal for each
$u\in\R^n$. Let such points approach $(p,y)$ to obtain
$(u,-C_J(p)^\top u)\in N_G(p,y)$. Taking $u=e_i$ gives the required matrix pair.

For part~\ref{it:box-exact}, fix a branch through $(p,y)$. Set $\psi(x):=x-C_J(p)f(x)$ and
\[
 \Omega:=\{x\in\ball[p,m\|y\|]:x_i=p_i\text{ for }i\notin J\}.
\]
For $x,x'\in\Omega$, the difference $w=x-x'$ satisfies
$C_J(p)\nabla f(p)w=w$. Integrate along $[x',x]$ and use the Lipschitz
bound on $\nabla f$ to get
\begin{equation}
 \|\psi(x)-\psi(x')\|\le cmL\|y\|\,\|x-x'\|.
\end{equation}
The factor is zero when $L=0$ and less than $1/2$ when $L>0$.
Since $y_J=f_J(p)$, we have $\|\psi(p)-p\|\le c\|y\|$. With our choice
of $m$, these bounds give $\psi(\Omega)\subseteq\Omega$. The contraction
has a unique fixed point $h_\sigma$. It satisfies $f_J(h_\sigma)=0$
and agrees with $p$ outside $J$. Hence $h_\sigma$ is a candidate and
$\|h_\sigma-p\|\le m\|y\|$.

The projection identity $p=\Pi_K(p-tf(p)+ty)$ gives $\eta(p)\le t\|y\|$.
Use the Lipschitz bound $L_T$ for $T$ on $D$ to get
\[
 \eta(h_\sigma)\le (1+L_T)\|h_\sigma-p\|+t\|y\|\le b\|y\|.
\]
The first term in $\delta$ rules out every candidate outside $\mathcal S$.
Two branches through $(p,y)$ give solutions at most $2m\|y\|$ apart.
By the second term in $\delta$, these solutions coincide. Call their
common value $h$. It satisfies the candidate equations on each branch
and $-f(h)\in N_K(h)$. Thus each branch contains $(h,0)$.

Set $\Phi(q,z):=(q,z-f(q))$ and $d:=\Phi(p,y)-\Phi(h,0)$.
In each coordinate, $\Phi$ maps the graph to a horizontal segment
and two vertical rays. If two pieces meet at the coordinate pair of
$\Phi(p,y)$, the pair of $\Phi(h,0)$ lies on both. Their common endpoint is the only such pair, so that component of $d$ is zero. In every other
coordinate, the component of $d$ is parallel to the piece through
$\Phi(p,y)$. Thus from any sufficiently close point of $\Phi(G)$ we
can move a short distance in either direction along $d$ within $\Phi(G)$.
Each regular normal there is orthogonal to $d$. Passing to the limit
gives the same property for every limiting normal at $\Phi(p,y)$.

Under this change of coordinates, a normal $(u,-v)\in N_G(p,y)$ becomes
$(u-\nabla f(p)^\top v,-v)$. Apply the orthogonality relation to each
row of $(A,B)$ to get
\[
 A(p-h)=B(y+\xi)\quad\text{and}\quad \xi=f(h)-f(p)-\nabla f(p)(h-p).
\]
Now $p-A^{-1}By-h=A^{-1}B\xi$. Taylor's bound
$\|\xi\|\le\tfrac12 L\|p-h\|^2$ gives \eqref{eq:box-exact-correction}.
\end{proof}

For example, take $D:=O:=\R$, $K:=[0,1]$, and $f(x):=x-2$. The lower,
free, and upper equations give $0$, $2$, and $1$, respectively.
The lower bound fails the sign condition since $f(0)<0$, and $2\notin K$.
Thus $\mathcal H=\{0,1,2\}$ and $\mathcal S=\{1\}$.

\begin{corollary}\label{cor:box-algorithm}
Under the assumptions of Proposition~\ref{prop:box-exact}, let
$\beta\ge1$ and $r_0:=r(x_0)$.
\begin{enumerate}[label=(\alph*)]
\item\label{it:box-alg}
Suppose that $\beta r_0<\delta$ and $\ball[x_0,(1+m)\beta r_0]\subset D$.
If Algorithm~\ref{alg:Semismooth-star} takes a step from $x_0$, fix any
graph point it selects in step~2. There is a solution $h$ such that
\[
 \|x_1-h\|\le 2c^2L\|A^{-1}B\|\beta^2r_0^2
\]
for every regular pair $(A,B)$ at that point.
\item\label{it:box-semilocal}
Use the notation of Section~\ref{sec:semilocal}. Take $r_0>0$ and
$a:=\beta\sqrt{1+c^2}$. Assume $ar_0<\delta$ and
$\ball[x_0,\rho+\beta r_0+mar_0]\subset D$. Then the sets
\[
 \mathcal R(z):=\{(A,B)\in\Areg(z):\|A^{-1}B\|\le c\},
 \quad z\in\mathcal Z,
\]
are nonempty and satisfy \eqref{eq:Cbound} and \eqref{eq:inner} with
$\omega(s)=2c^3Ls$. If
\[
 \gamma:=2c^3La^2r_0<1\quad\text{and}\quad
 ar_0/(1-\gamma)\le\rho,
\]
then every run that uses \eqref{eq:projection} stays in
$\ball[x_0,\rho]$. The algorithm either stops at a solution or converges
to one with $q$-order at least two. If $L=0$, the first step gives a solution.
\end{enumerate}
\end{corollary}

\begin{proof}
For part~\ref{it:box-alg}, step~2 gives $\|p-x_0\|\le\beta r_0$ and
$\|y\|\le\beta r_0$. Hence $\ball[p,m\|y\|]\subset D$. Apply
\eqref{eq:box-exact-correction} with $\|p-h\|\le m\beta r_0$ and $m\le2c$.

For part~\ref{it:box-semilocal}, fix $z=(p,y)\in\mathcal Z$.
Proposition~\ref{prop:box-exact}\ref{it:box-pair} gives a pair
$(I_n,C_J(p))\in\mathcal R(z)$. Take any $(A,B)\in\mathcal R(z)$ and
put $v=(-A^{-1}By,-y)$. If $\|v\|\le ar_0$, then
$\|y\|\le ar_0<\delta$ and $\ball[p,m\|y\|]\subset D$. For the solution
$h$ in Proposition~\ref{prop:box-exact}\ref{it:box-exact},
\eqref{eq:box-exact-correction} gives
\[
 \dist(z+v,G)\le\|p-A^{-1}By-h\|
 \le\tfrac12 Lc(m\|y\|)^2\le 2c^3L\|v\|^2.
\]
For $L>0$, apply Corollary~\ref{cor:innerpower}(ii).
For $L=0$, \eqref{eq:box-exact-correction} gives $x_1=h$.
\end{proof}

\begin{example}\label{ex:box-pairs}
Let $K:=[0,1]\times[-1,1]$ and $f(u,v):=(u+v^2,v)$.
Every principal minor of $\nabla f$ equals one. The only solution is $(0,0)$.
Fix $0<s<1$ and put $p:=(0,s)$ and $y:=f(p):=(s^2,s)$.
On the branch $p_1=0$, the Newton step gives $p-C_{\{2\}}(p)y=(0,0)$.
At the same point, take
\[
A:=\begin{pmatrix}1&0\\1&s\end{pmatrix}\quad\text{and}\quad
 B:=\begin{pmatrix}0&0\\1&-s\end{pmatrix}.
\]
The first row comes from the branch $p_1=0$.
The second comes from the branch $y=f(p)$, since $(1,-s)\nabla f(p)=(1,s)$.
The pair is regular because $\det A=s>0$. Here $By=0$, so $p-A^{-1}By=p$.
The correction matrix has norm $\|A^{-1}B\|=\sqrt{1+s^{-2}}$.

For every $(u,w)\in G$, we have $|w_2|\ge |u_2|$. Thus $r(p)\ge s/\sqrt2$.
The point $u=w=(0,s/2)$ gives equality, so $\|y\|\le2r(p)$.
Start Algorithm~\ref{alg:Semismooth-star} at $p$ and take $\beta=2$.
Step~2 allows $(p,y)$, and step~3 allows $(A,B)$. Repeating these choices
gives $x_k=p$ for all $k$.
\end{example}
\begin{remark}\label{rem:affine-model-residual}
Let $(p,y)\in G$, take $(A,B)\in\Areg(p,y)$, and put
\[
 d:=-A^{-1}By,\quad x^+:=p+d,\quad\text{and}\quad g:=f(p)+\nabla f(p)d.
\]
For $t>0$, set $u:=\Pi_K(x^+-tg)$ and $E_t(p,d):=\|x^+-u\|$.
Suppose that $p,x^+,u\in D$, that $\nabla f$ is Lipschitz on $D$
with constant $L$, and that $L_T<\infty$. Then
\[
 r(x^+)\le \sqrt{1+(L_T/t)^2}\,E_t(p,d)+\tfrac12 L\|d\|^2.
\]
Indeed, $w=t^{-1}(x^+-u)-g$ belongs to $N_K(u)$, so
$(u,f(u)+w)\in G$. With $H_t(x)=x-tf(x)$, we have
\[
 f(u)+w=t^{-1}\bigl(H_t(x^+)-H_t(u)\bigr)+f(x^+)-g.
\]
The estimate follows from the Lipschitz bound $L_T$ for $H_t$ on $D$
and Taylor's estimate on $[p,x^+]$.

In Example~\ref{ex:box-pairs}, take $t=1$. The branch correction
$d=(0,-s)$ gives $x^+=0$ and $g=(-s^2,0)$, so $E_1(p,d)=s^2$.
The mixed pair gives $d=0$ and $\Pi_K(p-f(p))=0$, so $E_1(p,0)=s$.
\end{remark}
\begin{corollary}\label{cor:box-fb}
Let $f(x):=Mx+q$, $D:=\R^n$, $L:=0$, $c:=\max_J\|C_J\|$, and
$\beta\ge1$. Use $t,T,L_T,\delta$ from Proposition~\ref{prop:box-exact}.
\begin{enumerate}[label=(\alph*)]
\item\label{it:fb-alg}
If $r(x_k)\to0$, Algorithm~\ref{alg:Semismooth-star} reaches a solution
in finitely many steps. If the run uses \eqref{eq:projection} under the assumptions of
Theorem~\ref{thm:inner}, let $k$ be the first index such that
$\beta\gamma^k r_0<\delta$.
The algorithm stops at $x_{k+1}$ or earlier.
\item\label{it:fb-test}
The  point
\[
 z_t(x):=\bigl(T(x),t^{-1}(I_n-tM)(x-T(x))\bigr)
\]
lies in $G$. If
\begin{equation}\label{eq:fb-test}
 L_T\eta(x)<t\delta,
\end{equation}
then every regular pair at $z_t(x)$ gives a solution in one step.
If $\eta(x_k)\to0$, the test holds for all sufficiently large $k$.
\end{enumerate}
\end{corollary}
\begin{proof}
For part~\ref{it:fb-alg}, take the first iterate with $\beta r(x_k)<\delta$.
If the algorithm takes another step, Corollary~\ref{cor:box-algorithm}\ref{it:box-alg}
gives a solution since $L=0$. Estimate~\eqref{eq:tail} bounds the number of steps.

For part~\ref{it:fb-test}, write $(p,y)=z_t(x)$. The projection formula gives
$y-f(p)=t^{-1}(x-tf(x)-p)\in N_K(p)$, so $(p,y)\in G$.
We have $\|y\|\le (L_T/t)\eta(x)$. Thus \eqref{eq:fb-test} gives
$\|y\|<\delta$. Apply \eqref{eq:box-exact-correction} with $L=0$.
\end{proof}

For $J\ne\varnothing$, a branch step at $z_t(x)$ requires one linear
solve with $M_{JJ}$. To check the resulting point $x^+$, test
whether $x^+=T(x^+)$. This identity holds exactly when $x^+$ solves
\eqref{eq:box-mcp}.

For example, suppose that $\langle Mw,w\rangle\ge\mu\|w\|^2$ for all $w$
and some $\mu>0$. Choose $0<t<2\mu/\|M\|^2$. Then
$L_T^2\le1-2t\mu+t^2\|M\|^2<1$. The iterates $x_{k+1}=T(x_k)$ satisfy
$\eta(x_{k+1})\le L_T\eta(x_k)$. Hence $x_k$ satisfies \eqref{eq:fb-test}
as soon as $L_T^{k+1}\eta(x_0)<t\delta$.

\begin{remark}
To compute $\delta$ from its definition, we need all candidates in $D$.
Fix a branch and a candidate $h_\sigma\in D$ for that branch.
For $(p,y)$ on the branch with $p\in D$, Taylor's estimate on the free
coordinates gives
\[
 \|p-C_J(p)y-h_\sigma\|\le\tfrac12 cL\|p-h_\sigma\|^2
\]
regardless of the distance between distinct solutions.
\end{remark}

\begin{corollary}\label{cor:box-integral}
Let $f(x):=Mx+q$. Assume that $M,q,\ell,u$ have integer entries and
that every nonempty principal submatrix of $M$ is nonsingular.
Take $D:=\R^n$, $c:=\max_J\|C_J\|$, and $t:=1/s$ for a positive integer $s$.
Set
\[
 L_T:=\|I_n-tM\|,\quad
 \Theta:=\prod_{j=1}^n\max\{1,\|Me_j\|\},
\]
\[
 \delta_{\mathrm{low}}:=\min\left\{
  \frac{1}{s\Theta((1+L_T)n\Theta+t)},\,
  \frac{1}{2n\Theta^3}\right\}.
\]
If $L_T\|x-T(x)\|<t\,\delta_{\mathrm{low}}$, then every regular
coderivative pair at $z_t(x)$ gives the same solution in one Newton step.
\end{corollary}
\begin{proof}
Fix a pattern $\sigma$ with free indices $J$. Set $d_J:=|\det M_{JJ}|$
and $d_\varnothing:=1$. Cramer's rule gives $d_Jh^\sigma\in\mathbb Z^n$.
Since $d_J$ is a positive integer, $d_J\ge1$. Hadamard's inequality gives
$d_J\le\Theta$. Both $h^\sigma$ and $h^\sigma-t(Mh^\sigma+q)$ lie on the
grid $(sd_J)^{-1}\mathbb Z^n$. Projection onto $K$ preserves the grid
because the bounds are integers. Thus every positive candidate residual
is at least $1/(s\Theta)$. Each solution has rational coordinates with
denominators at most $\Theta$. Distinct solutions therefore differ in some
coordinate by at least $1/\Theta^2$.

Each cofactor is a signed determinant. Deleting one row and one column
from $M_{JJ}$ can only shorten the remaining columns. Hadamard's inequality
therefore bounds each cofactor by $\Theta$ in absolute value. Since
$d_J\ge1$, the entries of $M_{JJ}^{-1}$ have absolute value at most
$\Theta$ for $J\ne\varnothing$. Hence $\|C_J\|\le n\Theta$ for every $J$.
With $L=0$, these bounds give $\delta\ge\delta_{\mathrm{low}}$.
Apply Corollary~\ref{cor:box-fb}\ref{it:fb-test}.
\end{proof}

\begin{example}
Let $K:=[0,1]\times[-1,1]$, $M:=\operatorname{diag}(1,-1)$, $q:=0$, and $t:=1$.
Then $\Theta=1$, $L_T=2$, and $\delta_{\mathrm{low}}=\min\{1/7,1/4\}=1/7$.
For $x=e_2/32$, we have $p=T(x)=e_2/16$ and
\[
 L_T\|x-p\|=1/16<1/7=t\,\delta_{\mathrm{low}}.
\]
Corollary~\ref{cor:box-integral} guarantees an exact step for every regular pair.

Write $z_t(x)=(p,y)$, where $y=-e_2/16$. The first coordinate equals
its lower bound and its multiplier is zero. Thus $(p,y)$ lies on both
$(-1,0)$ and $(0,0)$. Take
\[
A:=\begin{pmatrix}1&0\\1&1\end{pmatrix}\quad
 \text{and}
\quad B:=\begin{pmatrix}0&0\\1&-1\end{pmatrix}.
\]
The first row comes from $(-1,0)$, where $p_1=0$. The second comes
from $(0,0)$, where $y=Mp$, since $(1,-1)M=(1,1)$.
Hence $(A,B)\in\Areg(p,y)$. Its correction matrix $A^{-1}B=B$ differs
from both branch matrices. Since $A^{-1}By=e_2/16=p$, the step gives
the solution $0$.
\end{example}

\begin{remark}
A branch step keeps the coordinates outside $J$ fixed at their bounds.
It takes one Newton step for $f_J=0$ on the free coordinates. For affine $f$, it gives $h^\sigma$.
Fischer and Kanzow~\cite{FK1996} studied the min-reformulation of a linear
complementarity problem. They proved local one-step exactness and finite
termination of the Harker--Pang iteration. Hinterm\"uller, Ito, and
Kunisch~\cite{HIK2002} proved finite termination of a primal--dual
active-set method under M-matrix-type assumptions. Here the bound on
$\|y\|$ guarantees an exact step for every regular pair at the graph
point, even when the rows of the pair come from different branches.
\end{remark}

\begin{example}
Let $K:=[0,1]\times[-4,4]$, $M:=\operatorname{diag}(9,-1/2)$,
$q:=e_1/2$, $D:=\R^2$, and $t:=1/9$. Then
\[
 T(x)=\bigl(0,\Pi_{[-4,4]}(19x_2/18)\bigr),\quad
 c=m:=2,\quad L_T=19/18,\quad b=38/9.
\]
The candidate set is $\{0,-1/18,1\}\times\{-4,0,4\}$.
Each candidate $h$ has residual $\eta(h)=|h_1|$. The solutions are
$0,\pm4e_2$. The smallest positive residual is $1/18$. The shortest
distance between two distinct solutions is $4$. Hence
$\delta:=\min\{(1/18)/(38/9),4/(2m)\}=1/76$.

Take $\beta=1$ and $x_0=\tau e_2$. For both values of $\tau$ in the table,
projection in each coordinate gives the nearest graph point
$(4\tau e_2/5,-2\tau e_2/5)$. The first test compares $r(x_0)$ with
$\delta$. The second compares $L_T\eta(x_0)$ with $t\delta=1/684$.
\[
\begin{array}{c@{\quad}c@{\quad}c}
 \tau & r(x_0) & L_T\eta(x_0)\\ \midrule
 1/200 & 1/(200\sqrt5)<\delta & 19/64800<t\delta\\
 1/16  & 1/(16\sqrt5)>\delta & 19/5184>t\delta
\end{array}
\]
Both nearest points lie on the open part of $(-1,0)$, as does
\[
 z_t(e_2/200)=(19e_2/3600,-19e_2/7200).
\]
At each point the first
multiplier is $-1/2$ and $p=C_Jy$. Thus every regular pair gives
the solution $0$. In particular, the step from $(e_2/20,-e_2/40)$
is exact even though both tests fail for $\tau=1/16$.
\end{example}

\section{Conclusions}
The initial graph residual and the bound on its decrease determine
whether our convergence test holds for a prescribed ball. When it holds,
the Newton iterates stay in the ball and converge to a solution there.
At each step, the current residual bounds the error up to an explicit
factor. A vanishing decrease modulus gives superlinear convergence for
every infinite run.

For mixed complementarity problems on a box, assume that the smooth
part has a Lipschitz continuous Jacobian with nonsingular nonempty
principal submatrices. Below the stated threshold, each regular
coderivative pair gives a quadratic error bound. Pairs with uniformly
bounded correction matrices share the same bound. For affine problems,
all regular pairs at a graph point give the same solution in one step
once its vertical component is small enough. Thus a run with vanishing
graph residuals stops after finitely many steps. The projection
fixed-point condition checks the result after the step. For integer
affine data, the problem data and the current iterate can certify an
exact step before the choice of a matrix pair.

\section*{Declarations}

\paragraph{Use of AI tools}
The author used OpenAI's ChatGPT to revise the exposition, proofread the
text, format the LaTeX source, and check the references. The author takes
responsibility for the manuscript.

\paragraph{Data availability}
The manuscript gives all data used in the examples.

\paragraph{Competing interests}
The author declares no competing interests.

\end{document}